# Physical interpretation of the Riemann hypothesis

Dmitry Pozdnyakov

*Faculty of Radiophysics and Computer Technologies of Belarusian State University,
Nezavisimosty av.4, Minsk 220030, Belarus
E-mail: pozdnyakov@tut.by*

**Abstract:** An equivalent formulation of the Riemann hypothesis is given. The physical interpretation of the Riemann hypothesis equivalent formulation is given in the framework of quantum theory terminology. One more power series related to the Riemann Xi function and the Riemann hypothesis is considered. Some roots of the polynomial connected with the power series are studied. It is shown that the Riemann hypothesis is true. But it is undecidable and must be considered as an axiom.

As the Riemann hypothesis so the problem of its proof are so well-known that writing of even a short introduction is unreasonable rewriting of copy-book maxims. But it is only necessary to note that a new point of view on the hypothesis in the framework of physical applications of the Riemann zeta function (see, for example, [1–3]) is proposed in this study.

**Mathematics**

Like Riemann let us exclude from consideration the trivial zeroes of the Riemann zeta function $\zeta(s) = \zeta(\sigma + it)$ for convenience. Their existence is caused by the fact that the Euler $\Gamma$-function entering the expression for $\zeta$-function [4] has singular points. Therefore instead of the Riemann zeta function let us consider Riemann's upper-case Xi function

$$\Xi(t,\sigma) = \Xi(\sigma + it) = \Xi(s) = \pi^{-s/2}\Gamma(s/2)\zeta(s) = \frac{1}{s(s-1)} + \int_1^\infty \left( x^{\frac{s}{2}-1} + x^{\frac{1-s}{2}-1} \right) \theta(x)dx , \quad (1)$$

where

$$\theta(x) = \sum_{n=1}^\infty \exp(-\pi n^2 x)$$

is the theta series (theta function). The zeroes of Riemann $\Xi$-function evidently coincide with the nontrivial zeroes of Riemann $\zeta$-function [1, 4]. Let us introduce new variables $\delta = \sigma - 1/2$ and $y = \ln(x)$. After a number of transformations of eq.(1) one can obtain relation

$$\Xi(t,\delta) = \Xi(t - i\delta) = 2\int_0^\infty \cos\left(\frac{t-i\delta}{2}y\right)\Omega(y)dy - \frac{1}{(t-i\delta - i/2)(t-i\delta + i/2)} =$$

$$2\int_0^\infty \cos(zy/2)\Omega(y)dy - \frac{1}{z^2 + 1/4} = \Xi(z) , \quad (2)$$

where

$$\Omega(y) = \sum_{n=1}^\infty \exp\left((y/4) - \pi n^2 \exp(y)\right).$$

So far as it is well known that the zeta function has not the nontrivial zeroes outside the critical strip $|\delta| > 1/2$ and on its boundaries $|\delta| = 1/2$ (J. Hadamard, Ch. J. La Vallee Poussin, 1896) let us consider only the case when values of $\delta$ are in the interval from $-1/2$ to $1/2$ ($|\delta| < 1/2$). Let us apply relation

$$\int_0^\infty \cos\left(\frac{t-i\delta}{2}y\right)\exp\left(-\frac{y}{4}\right)dy = \frac{1}{(t-i\delta - i/2)(t-i\delta + i/2)} \quad (3)$$

for subsequent transformation of eq.(2). It is true $\forall t \in \mathbb{R}$ and $\forall \delta \in (-1/2, 1/2)$. Eq.(3) is the result of Fourier transform of the second term in eq.(2). Then we get expression

$$\Xi(z) = 2\int_0^\infty \cos(zy/2)\Phi(y)dy , \quad (4)$$

where
$$\Phi(y) = -\frac{1}{2}\left[\exp(-y/4) - 2\sum_{n=1}^{\infty}\exp\left((y/4) - \pi n^2 \exp(y)\right)\right].$$

It should be noted that identical equality
$$\int_{-\infty}^{+\infty} g(y)\exp(\pm izy/2)\,dy \equiv 2\int_0^{\infty} g(y)\cos(zy/2)\,dy \quad (5)$$

is true in the sense of Cauchy principal value of the integrals for every even function $g$. Taking into account the limits of integration in eq.(4) the function $\Phi(y)$ can be naturally replaced by the function $\Phi(|y|)$ for which the identity is true. As a result we have equation
$$\Xi(z) = \int_{-\infty}^{+\infty} \exp(\pm izy/2)\Phi(|y|)\,dy \quad (6)$$

instead of eq.(4). The choice of sign before the imaginary unit does not obviously influence the zeroes of $\Xi$-function (see eq.(5)). Let us choose the lower sign for definiteness.

Let us change the variables $2t \to k$, $2\delta \to \lambda$, $y/4 \to x$ and introduce function
$$\rho_R(k,\lambda) = |A_R(k,\lambda)|^2 = \left|-\frac{\Xi(k,\lambda)}{2\sqrt{2\pi}}\right|^2 = \frac{1}{2\pi}\left|\int_{-\infty}^{+\infty} R(x)\exp(-i(k-i\lambda)x)\,dx\right|^2 =$$
$$\frac{1}{2\pi}\left|\int_{-\infty}^{+\infty} R(x)\exp(-iKx)\,dx\right|^2 = \left|-\frac{\Xi(K)}{2\sqrt{2\pi}}\right|^2 = |A_R(K)|^2 = \rho_R(K), \quad (7)$$

where
$$R(x) = \begin{cases} r^-(x), & x < 0; \\ \left(r^-(x) + r^+(x)\right)/2, & x = 0; \\ r^+(x), & x > 0; \end{cases} \quad (8)$$

$$r^-(x) = \exp(x) - 2\sum_{n=1}^{\infty}\exp\left(-x - \pi n^2 \exp(-4x)\right),$$

$$r^+(x) = \exp(-x) - 2\sum_{n=1}^{\infty}\exp\left(x - \pi n^2 \exp(4x)\right).$$

The zeroes of $\rho_R$-function $\{K_n = (k-i\lambda)_n \mid n \in \mathbb{N}\}$ obviously coincide with the zeroes of $\Xi$-function $\{z_n = (t-i\delta)_n \mid n \in \mathbb{N}\}$ up to constant factor, that is $K_n = 2z_n$, $n \in \mathbb{N} = \{1,2,3,...\}$. But, in contrast to $\Xi$-function, all the values of $\rho_R$-function are obviously the real positive numbers or zero, in other words $\forall K \in \mathbb{C}\ \rho_R \geq 0$.

So, basing on eqs.(1) – (8), the Riemann hypothesis is reduced to (completely equivalent to) conjecture that $\forall k \in \mathbb{R}$ and $\forall \lambda \in (-1,1) \setminus \{0\}\ \rho_R > 0$. Thus the Riemann hypothesis given in such a form is a stronger assertion than the evident non-strict inequality $\rho_R \geq 0\ \forall k \in \mathbb{R}$ and $\forall \lambda \in (-1,1)$. If the Riemann hypothesis was false it would mean that $\exists k \in \mathbb{R}$ and $\exists \lambda \in (-1,1) \setminus \{0\}: \rho_R = 0$.

**Physics**

We need some generalizations for the physical interpretation of the hypothesis. Let us consider a set of wave functions $\{\psi\}$ in the coordinate representation $\psi: \mathbb{R} \to \mathbb{C}$ which are the stationary Schrödinger equation solutions describing some bounded quantum states [5, 6]. Then by analogy with eq.(7) we have relation
$$\rho(k,\lambda) = |A(k,\lambda)|^2 = \frac{1}{2\pi}\left|\int_{-\infty}^{+\infty} \psi(x)\exp(-i(k-i\lambda)x)\,dx\right|^2 =$$
$$\frac{1}{2\pi}\left|\int_{-\infty}^{+\infty} \psi(x)\exp(-iKx)\,dx\right|^2 = |A(K)|^2 = \rho(K). \quad (9)$$

It is evident that $\forall K \in \mathbb{C}\ \rho \geq 0$. In context of quantum physics terminology in eq.(9) $K$ is the wave vector which is complex in general case ($K = k - i\lambda$; $k, \lambda \in \mathbb{R}$) [5, 6]; $x$ is the coordinate (position); $A$ is the wave function in the wave vector representation [5, 6]; $\rho$ is the function of spectral density of quantum states [5, 6]. The considered functions are wave functions describing quantum states of particles in case of their nonrelativistic one-dimensional finite motion.

If $\lambda = 0$ then eq.(9) characterizes the expansion of a $\psi$-function in the stable states of wave-vector space $\mathbb{R}_k$ which is a subspace of phase space $\mathbb{R}_x \times \mathbb{R}_k$. Here $\mathbb{R}_x$ is the coordinate (position) space which is a subspace of phase space too. The formalism of path integrals can be applied to describe any quantum system with the stable (non-decaying) states [7]. It is equivalent to formalism based on an evolution wave equation with a real Hamiltonian [7]. The quantities $k_n$ ($n \in \mathbb{N}$) satisfying the equality $\rho(k_n) = 0$ are nothing else than forbidden (non-excited) states in the wave-vector space. If $\lambda \neq 0$ then eq.(9) characterizes the expansion of a $\psi$-function in the unstable (decaying) states of complex wave-vector space $\mathbb{C}_K$ which is a subspace of generalized phase space $\mathbb{R}_x \times \mathbb{C}_K$. As a rule an unstable quantum state is either the spontaneously decaying quantum state or the quantum state decaying during irreversible decoherence process. The formalism of restricted path integrals can be applied to describe any quantum system with the unstable (decaying) states [7]. Under some conditions it is equivalent to formalism based on an evolution wave equation with a complex Hamiltonian [7]. Thus, during passage from $k \in \mathbb{R}$ ($\lambda = 0$) to $K \in \mathbb{C} \setminus \mathbb{R}$ ($\lambda \neq 0$) we are passing on from consideration of physical systems with the stable quantum states to consideration of physical systems with the unstable ones. In particular, the passage is equivalent to passage from consideration of isolated quantum systems to consideration of open ones [7].

It is easy to show (see Appendix A) that the function $R(x)$ can be considered as the wave function. Let us call it as the Riemann wave function (the Riemann wave function in the coordinate representation). Let us also call the function $A_R(K)$ as the Riemann wave packet (the Riemann wave function in the wave vector representation), and the function $\rho_R(K) = |A_R(K)|^2$ as the Riemann spectral function (function of spectral density of quantum states).

Let us now appeal to a simple but significant example of finite motion of a particle that will sufficiently illustrate the disappearance of forbidden states in the spectrum at passage from the expansion of $\psi$-functions in the real wave vectors $k$ to the expansion of $\psi$-functions in the complex wave vectors $K = k - i\lambda$. In particular, let us consider a particle in an infinitely deep potential well of width $a$ for which the states are defined by the wave functions [5, 6]

$$\psi_n(x) = \begin{cases} \sqrt{2/a} \cos(n\pi x/a), x \in [-a/2, a/2], n = 1, 3, \ldots; \\ \sqrt{2/a} \sin(n\pi x/a), x \in [-a/2, a/2], n = 2, 4, \ldots; \\ 0, x \notin [-a/2, a/2], n = 1, 2, 3, 4, \ldots. \end{cases}$$

In that case the spectral density calculated by means of eq.(9) is given by expression

$$\rho_n(k, \lambda) = \begin{cases} \dfrac{4\pi n^2}{a^3} \dfrac{\cos^2(ka/2)\cosh^2(\lambda a/2) + \sin^2(ka/2)\sinh^2(\lambda a/2)}{\left(k^2 - \lambda^2 - \pi^2 n^2 a^{-2}\right)^2 + 4\lambda^2 k^2}, \\ k \in \mathbb{R}, \lambda \neq 0, n = 1, 3, \ldots; \\ \dfrac{4\pi n^2}{a^3} \dfrac{\sin^2(ka/2)\cosh^2(\lambda a/2) + \cos^2(ka/2)\sinh^2(\lambda a/2)}{\left(k^2 - \lambda^2 - \pi^2 n^2 a^{-2}\right)^2 + 4\lambda^2 k^2}, \\ k \in \mathbb{R}, \lambda \neq 0, n = 2, 4, \ldots; \\ \dfrac{4\pi n^2}{a^3} \dfrac{\cos^2(ka/2)}{\left(k^2 - \pi^2 n^2 a^{-2}\right)^2}, k \neq \pm\pi n a^{-1}, \lambda = 0, n = 1, 3, \ldots; \\ \dfrac{4\pi n^2}{a^3} \dfrac{\sin^2(ka/2)}{\left(k^2 - \pi^2 n^2 a^{-2}\right)^2}, k \neq \pm\pi n a^{-1}, \lambda = 0, n = 2, 4, \ldots; \\ \dfrac{a}{4\pi}, k = \pm\pi n a^{-1}, \lambda = 0, n = 1, 2, 3, 4, \ldots. \end{cases}$$

It is evident that $\forall k \in \mathbb{R}$, $\forall \lambda \in \mathbb{R} \setminus \{0\}$, $\forall n \in \mathbb{N}$ $\rho_n > 0$, and $\forall m, n \in \mathbb{N}$ $\rho_n(\pm \pi m a^{-1}, 0) = 0$ at $m \neq n$. The function $\rho_n(k, 0)$ vanishes in an infinite number of points like the function $\rho_R(k, 0)$.

So, taking into account everything mentioned above, an equivalent formulation of the Riemann hypothesis in the framework of quantum theory terminology can be laid down, for example, like this: *the Riemann spectral function vanishes only for the real values of wave vector*.

**Some more mathematics**

Expanding the cosine in eq.(2) in the Taylor series, the first term in the equality can be represented by a power series. That is, eq.(2) has turned into the following equality

$$\Xi(z) = 2\left[\sum_{m=0}^{\infty} c_m z^{2m}\right] - \frac{1}{z^2 + 1/4}, \quad (10)$$

where

$$c_m = \frac{(-1)^m}{4^m (2m)!} \sum_{n=1}^{\infty} \int_0^{\infty} y^{2m} \exp\left((y/4) - \pi n^2 \exp(y)\right) dy. \quad (11)$$

According to ref.[8] the improper integral in eq.(11) is reduced to the corresponding order derivative of the incomplete Gamma-function. Namely, we have the expression

$$\int_0^{\infty} y^{2m} \exp\left((y/4) - \pi n^2 \exp(y)\right) dy = \int_1^{\infty} t^{-3/4} \ln^{2m}(t) \exp\left(-\pi n^2 t\right) dt =$$

$$\frac{d^{2m}}{dx^{2m}} \left[(\pi n^2)^{-x} \Gamma(x, \pi n^2)\right]\bigg|_{x=1/4}.$$

And as a result we come to the formula

$$c_m = \frac{(-1)^m}{4^m (2m)!} \sum_{n=1}^{\infty} \frac{d^{2m}}{dx^{2m}} \left[\exp(-x \ln(\pi n^2)) \Gamma(x, \pi n^2)\right]\bigg|_{x=1/4}. \quad (12)$$

According to Cauchy–Hadamard theorem the series in eq.(10) converges $\forall z \in \mathbb{C}$. In particular, anyone can be convinced that $\forall m \in \mathbb{N}$ $|c_m| < \left(e \ln(1 + 2m/\pi)/(4m)\right)^{2m} \Rightarrow |c_m|^{1/(2m)} \to 0$ at $m \to \infty$.

It is evident (see eq.(10)) that zeros of $\Xi(z)$ completely coincide with zeros of the power series

$$P_{\infty}(z) = \frac{1}{2} - (z^2 + 1/4) \sum_{m=0}^{\infty} \frac{(-1)^m a_m z^{2m}}{4^m (2m)!}, \quad (13)$$

where

$$a_m = \sum_{n=1}^{\infty} \int_0^{\infty} \exp\left((y/4) + 2m \ln(y) - \pi n^2 \exp(y)\right) dy = \sum_{n=1}^{\infty} \frac{d^{2m}}{dx^{2m}} \left[\exp(-x \ln(\pi n^2)) \Gamma(x, \pi n^2)\right]\bigg|_{x=1/4}. \quad (14)$$

Taking into account the explicit form of eqs.(2), (10) and (13), the Riemann hypothesis is also equivalent to one more assumption: *all the roots of equation*

$$P_{\infty}(z) = 0 \quad (15)$$

*are the real numbers*.

Let us further consider the following polynomial

$$P_M(z) = \frac{1}{2} - (z^2 + 1/4) \sum_{m=0}^{M} \frac{(-1)^m a_m z^{2m}}{4^m (2m)!}, \quad (16)$$

where $M = 0, 1, 2, \ldots$, which relates to the power series $P_{\infty}(z)$. For large values of $M$, as is known [9], zeroes of $P_M(z)$ can be found only by means of numerical methods.

In figs.1 and 2 the roots $z_n$ of the polynomial equation

$$P_M(z) = 0, \quad (17)$$

are represented for some values of $M$.

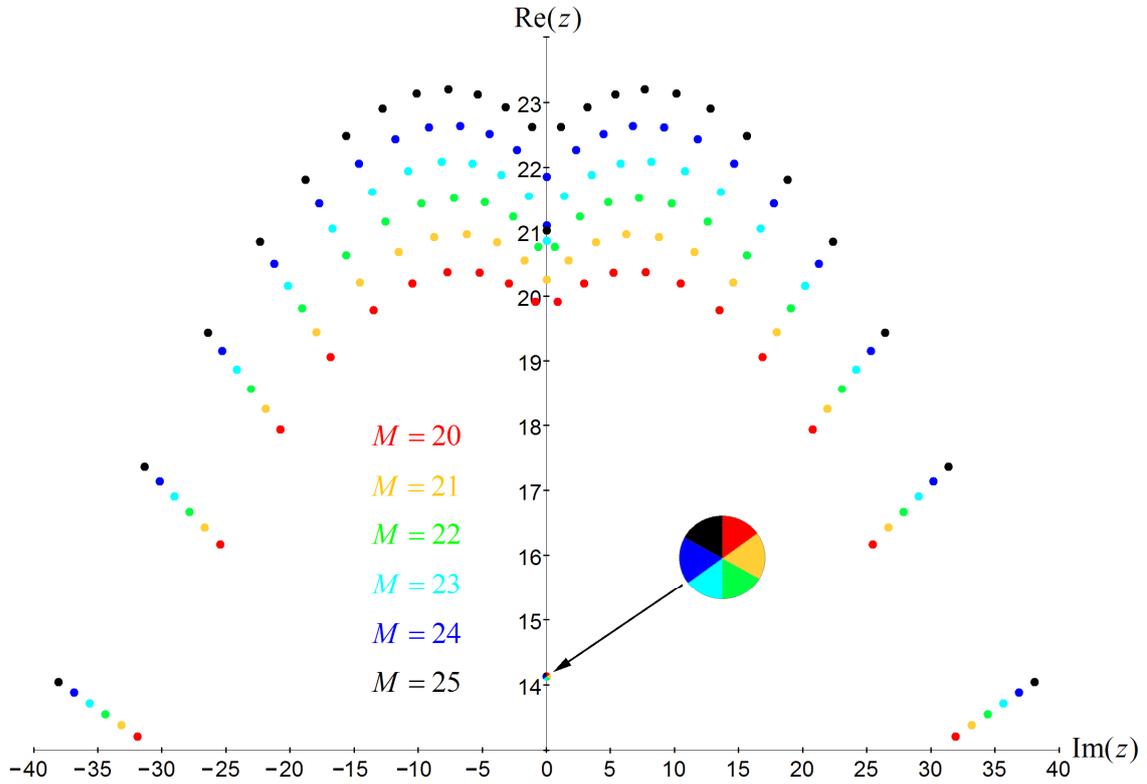

**Fig.1** All roots of eq.(17) in the complex half-plane $\mathrm{Re}(z) > 0$.
The numerical error in finding the roots is less than 0.01, i.e. $\forall z_n \quad |z_n^{\mathrm{numerical}} - z_n^{\mathrm{exact}}| < 0.01$.

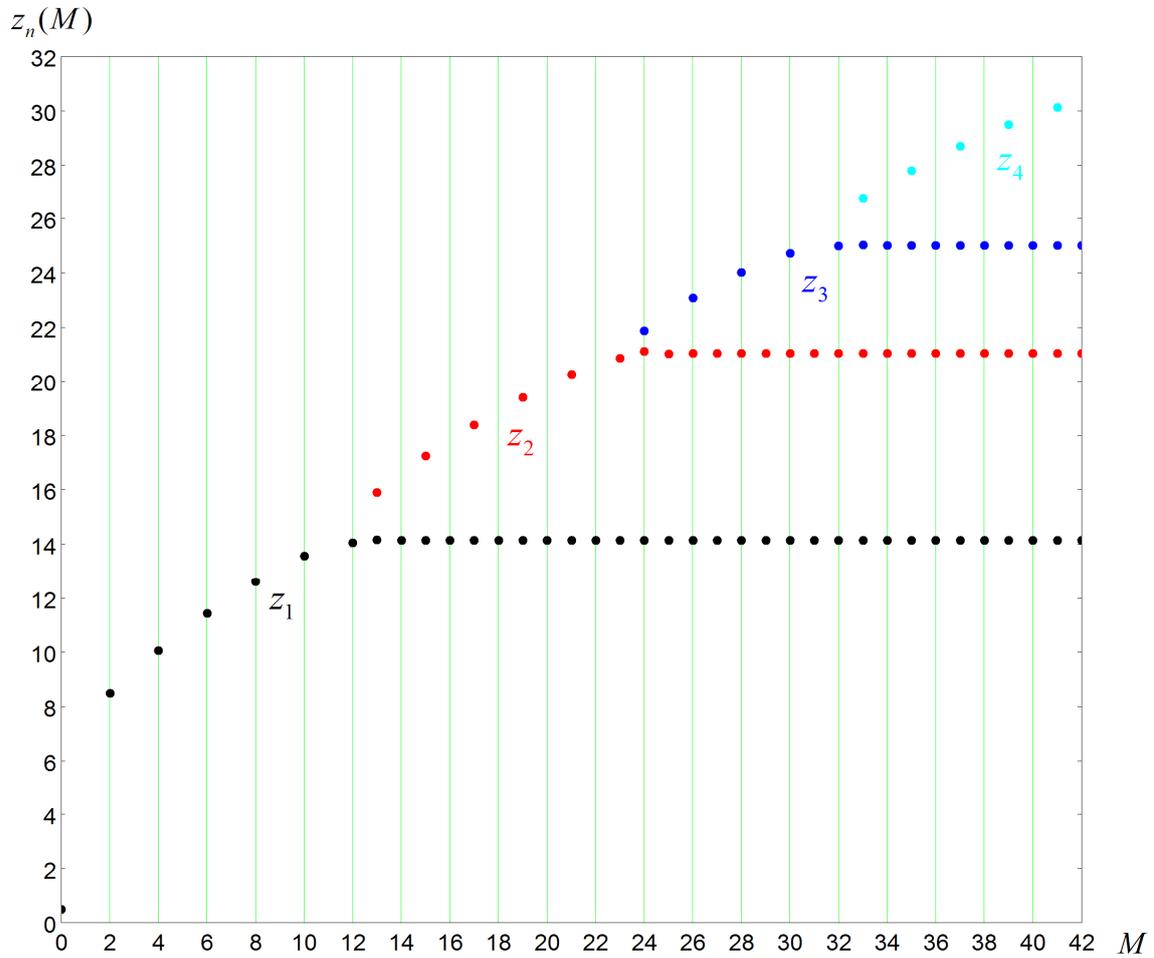

**Fig.2** All real positive roots of eq.(17).
The numerical error in finding the roots is less than $10^{-4}$, i.e. $\forall z_n \quad |z_n^{\mathrm{numerical}} - z_n^{\mathrm{exact}}| < 10^{-4}$.

It is evident from the figures that in contrast to the real roots of eq.(17) there is no disorder in distribution of the complex roots. They are strictly ordered.

**Important remarks**

In accordance with ref. [4] the Riemann Xi function $\Xi(s)$, and therefore the function $\Xi(z)$, is not represented by an analytical expression in closed form, i.e. it is no way to pass from the integral in eq.(2) to an expression with a finite number of operations over a finite number of elementary functions (the integral is not represented by quadratures). As a result there is no analytical way to find and verify the roots of eq.(2). Also there is no analytical way to find and verify the roots of eq.(15) because of both the impossibility to represent coefficients $a_m$ by quadratures and the infinite number of terms in eq.(15) [9]. In that case it is impossible to check analytically the belonging of every root $z_n$ to $\mathbb{R}$ or $\mathbb{C} \setminus \mathbb{R}$ ($\forall z_n \in \mathbb{C}$). Finally the only way of such a search and examination is a numerical one. But, because of the infinite number of roots of eqs.(2) and (15) [4], any numerical search and examination of them all is impossible from the practical point of view since it takes infinitely long computational time. Consequently, we have two situations. The first one corresponds to the case when the Riemann hypothesis is true. Then obviously every found root is a real number. At that the rest part of the roots, which is always infinite, has never been checked. It follows from this that if the hypothesis is true then it is impossible to prove it. Naturally the hypothesis negation cannot be proved for such a case in principle. The second situation corresponds to the case when the Riemann hypothesis is false, and consequently it is unprovable in principle, but its negation is provable since there is an algorithm which will find at least one root $z^* \in \mathbb{C} \setminus \mathbb{R}$ sooner or later. On this basis we can conclude that the Riemann hypothesis can only be either consistent and undecidable or decidable and inconsistent.

According to ref. [10] the Riemann hypothesis is consistent and therefore it is undecidable. As a result it is true according to Gödel's incompleteness theorems. For such a case we could draw a parallel with the Euclidean geometry (see Euclid's fifth postulate: the parallel postulate) when considering the Riemann hypothesis as "the fifth postulate of mathematics" *sui generis*. In particular, Gödel's incompleteness of *a priori* incomplete axioms system, in which the Riemann hypothesis is formulated, could appear itself just through the hypothesis. At that the system could be completed if to consider the hypothesis as the statement, i.e. missing postulate or axiom, and not the assumption.

**Conclusions**

Thus, the equivalent formulation of the Riemann hypothesis is given. The physical interpretation of the Riemann hypothesis equivalent formulation is given in the framework of quantum theory terminology. The power series $P_\infty(z)$ related to the Riemann Xi function and the Riemann hypothesis is considered. Some roots of the polynomial $P_M(z)$ connected with the power series $P_\infty(z)$ are studied. It is shown that the Riemann hypothesis is undecidable and true according to Gödel's incompleteness theorems.

**Acknowledgement**

The author is very grateful to Dr. Serguei K. Sekatski (Laboratory of Physics of Living Matter, Ecole Polytechnique Fédérale de Lausanne, Lausanne, Switzerland) for helpful discussions and a critical review of the paper.

**Appendix A**

$R(x)$ (see eq.(8)) is the wave function as it satisfies the basic requirements of quantum mechanics (the regularity conditions) [5, 6], and it is the eigenfunction of the Schrödinger Hamiltonian. Let us show this below.

1. $R(x)$ is evidently single-valued.
2. $R(x)$ is evidently finite.
3. $R(x)$ is normalizable since

$$\int_{-\infty}^{+\infty} |R(x)|^2 \, dx < \int_{-\infty}^{+\infty} |\exp(-|x|)|^2 \, dx = 1.$$

4. The function $R(x)$ and its first derivative $R'(x)$ are obviously continuous $\forall x \in \mathbb{R} \setminus \{0\}$ (see eq.(8)). $R(x)$ and $R'(x)$ are also continuous in the point $x = 0$ since

$$\lim_{x \to 0} R(x) = 1 - 2\sum_{n=1}^{\infty} \exp(-\pi n^2) = R(0),$$

$$\lim_{x \to 0} R'(x) = 0 = R'(0).$$

Anyone can be convinced of this fact, applying the well-known relation for the theta series $\theta(x)$ [4]
$$4\theta'(1) + \theta(1) = -1/2.$$

5. It is easy to be convinced that the function $R(x)$ satisfies the Schrödinger equation
$$-R''(x) + u_R(x)R(x) = \varepsilon_R R(x),$$
where
$$\varepsilon_R = -1,$$
$$u_R(x) = 16\pi \frac{\sum_{n=1}^{\infty} n^2 \left(3 - 2\pi n^2 \exp(4|x|)\right) \exp\left(6|x| - \pi n^2 \exp(4|x|)\right)}{1 - 2\sum_{n=1}^{\infty} \exp\left(2|x| - \pi n^2 \exp(4|x|)\right)}.$$

Thus $R(x)$ is the wave function describing a bound state. Namely, the ground quantum state in the potential well characterized by the potential function $u_R$.

**Appendix B**

It is extremely important and interesting to search and investigate stochastic processes of a certain kind which apparently have not been studied yet. The processes are non-Markovian stochastic processes described by the following autocorrelation function
$$\tau(t) = \frac{1 - 2\sum_{n=1}^{\infty} \exp\left(2|t|t_0^{-1} - \pi n^2 \exp(4|t|t_0^{-1})\right)}{1 - 2\sum_{n=1}^{\infty} \exp(-\pi n^2)} \exp\left(-|t|t_0^{-1}\right),$$
where $t_0 \in \mathbb{R}_+$.